\documentclass[number]{elsarticle}
\usepackage[english]{babel}
\usepackage{amsmath, amsthm}
\usepackage{amssymb}
\usepackage{latexsym}
\usepackage{fancyhdr}

\newcommand{\K}{\mathcal{K}}
\newcommand{\F}{\mathcal{F}}
\newcommand{\D}{\mathcal{D}}
\newcommand{\C}{\mathcal{C}}

\newcommand\PG{\mathrm{PG}}
\newcommand\PGL{\mathrm{P\Gamma L}}

\newcommand\GF{\mathrm{GF}}
\newcommand\Aut{\mathrm{Aut}}
\newcommand\ggd{\mathrm{ggd}}
\newcommand\id{\mathrm{id}}

\newcommand\Tr{\mathrm{Tr}}
\newcommand\V{\mathrm{V}}
\newcommand\N{\mathbb{N}}

\def\bew{{\bf Proof. }\ }
\def\evb{\hfill \ensuremath{\Box}\\ }

\newtheorem{lemma}{Lemma}
\newtheorem{theorem}{Theorem}
\newtheorem{remark}{Remark}
\newtheorem{corollary}{Corollary}

\title{A geometric approach to Mathon maximal arcs}
\date{}

\author{F. De Clerck\corref{cor1}}
\ead{fdc@cage.ugent.be} 
\address{ Department of Pure Mathematics and Computer Algebra, 
Ghent University, Krijgslaan 281 - S22, B-9000 Gent, BELGIUM}

\author{S. De Winter\fnref{fn1}} 
\ead{sgdwinte@cage.ugent.be} 
\address{ Department of Mathematics, Ohio University, Morton 321, Athens OH 45701, USA}

\author{T. Maes} 
\ead{tmmaes@cage.ugent.be} 
\address{ Department of Pure Mathematics and Computer Algebra, 
Ghent University, Krijgslaan 281 - S22, B-9000 Gent, BELGIUM}

\cortext[cor1]{Corresponding author} 
\fntext[fn1]{Post-doctoral Fellow of the Science Foundation Flanders(FWO-Vlaanderen)} 

\begin{document}
\begin{frontmatter}

\begin{abstract}
In 1969  Denniston gave a construction of maximal arcs of degree $d$ in Desarguesian projective planes of even order $q$, for all $d$ dividing $q$.  In 2002 Mathon gave a construction method generalizing the one of Denniston.  We will give a new geometric approach to these maximal arcs.  This will  allow us to count the number of isomorphism classes of Mathon maximal arcs of degree $8$ in $\PG(2,2^h)$, $h$ prime.
\end{abstract}
\end{frontmatter}

\section{Introduction}\label{intro}

A $\{k;d\}$-arc $\K$ in a finite projective plane of order $q$ is a non-empty proper subset of $k$ points such that some line of the plane meets $\K$ in $d$ points, but no line meets $\K$ in more than $d$ points.  For given $q$ and $d$, $k$ can never exceed $q(d-1)+d$.  If equality holds $\K$ is called a \textit{maximal arc of degree $d$}, \emph{a degree $d$ maximal arc}, \emph{a $\{q(d-1)+d;d\}$-arc} or shorter, \emph{a $d$-arc}.  Equivalently, a maximal arc can be defined as a non-empty, proper subset of points  of a projective plane, such that every line meets the set in $0$ or $d$ points, for some $d$.  The set of points of an affine subplane of order $d$ of a projective plane of order $d$ is a trivial example of a $\{d^2;d\}$-arc, as well as a single point, being a $\{1;1\}$-arc of the  projective plane.  We will neglect for the rest of this paper these trivial examples. \par
If $\K$ is a $\{q(d-1)+d;d\}$-arc in a projective plane of order $q$, the set of lines external to $\K$ is a $\{q(q-d+1)/d;q/d\}$-arc in the dual plane.  It follows that a necessary condition for the existence of a $\{q(d-1)+d;d\}$-arc in a projective plane of order $q$ is that $d$ divides $q$.  Denniston \cite{MR0239991} showed that this necessary condition is sufficient in the Desarguesian projective plane $\PG(2,q)$ of order $q$ when $q$ is even.  Ball, Blokhuis and Mazzocca \cite{MR1466573} showed that no non-trivial maximal arcs exist in a Desarguesian projective plane of odd order.  Note that if $\pi$ is a Desarguesian plane of order $q$ which contains a maximal arc $\K$ of degree $d$, then it also contains a maximal arc of degree $q/d$, the so-called \emph{dual maximal arc} of $\K$. \par
In \cite{MR1883870}, Mathon constructed maximal arcs in Desarguesian projective planes generalizing the previously known construction of Denniston \cite{MR0239991}.  We will begin by describing this construction of Mathon.  

From now on let $q=2^h$ and let $\Tr$ denote the usual absolute trace map from the finite field $\GF(q)$ onto $\GF(2)$.  We represent the points of the Desarguesian projective plane $\PG(2,q)$ as triples $(a,b,c)$ over $\GF(q)$, and the lines as triples $[u,v,w]$ over $\GF(q)$. A point $(a,b,c)$ is incident with a line $\left[u,v,w\right]$ if and only if $au+bv+cw=0$.  For $\alpha, \beta \in \GF(q)$ such that $\Tr(\alpha \beta)=1$, and $\lambda \in \GF(q)$ we define 
\[
F_{\alpha,\beta,\lambda}=\{(x,y,z):\alpha x^2+xy+\beta y^2+\lambda z^2=0\}.
\]  
Remark that $F_{\alpha,\beta,\lambda}$ is a conic if $\lambda \neq 0$, and that all the conics  have the point $F_{\alpha,\beta,0}=F_0=(0,0,1)$ as their nucleus. Due to the trace condition, the line $z=0$ is external to all conics.  \par
Let $\F$ be the set of all  $F_{\alpha,\beta,\lambda}$,   $\lambda \in \GF(q)$. For given $\lambda \neq \lambda'$, define a composition
\[
F_{\alpha, \beta, \lambda} \oplus F_{\alpha', \beta', \lambda'}=F_{\alpha \oplus \alpha', \beta \oplus \beta', \lambda \oplus \lambda'}
\]
where the operator $\oplus$ is defined  as follows:
\[
\alpha \oplus \alpha'= \frac{\alpha \lambda + \alpha'\lambda'}{\lambda + \lambda'}, \>\> \beta \oplus \beta'= \frac{\beta\lambda + \beta'\lambda'}{\lambda + \lambda'}, \>\> \lambda \oplus \lambda'= \lambda + \lambda'.
\]
The following lemma was proved by Mathon in \cite{MR1883870}.

\begin{lemma}  \label{discon}
Two non-degenerate conics $F_{\alpha,\beta,\lambda}$, $F_{\alpha',\beta',\lambda'}$, $\lambda \neq \lambda'$ and their composition $F_{\alpha,\beta,\lambda} \oplus F_{\alpha',\beta',\lambda'}$ are mutually disjoint if $\Tr((\alpha \oplus \alpha')(\beta \oplus \beta'))=1$.
\end{lemma}

 Given some subset $\C$ of $\F$, we say $\C$ is \textit{closed} if for every $F_{\alpha,\beta,\lambda} \neq F_{\alpha',\beta',\lambda'} \in \C$, $F_{\alpha \oplus \alpha',\beta \oplus \beta',\lambda \oplus\lambda'} \in \C$.  We can now state Mathon's theorem.

\begin{theorem} [\cite{MR1883870}]  \label{mathon}
Let $\C\subset\F$ be a closed set of conics  in $\PG(2,q)$, $q$ even.  Then the union of the points on the conics of $\C$ together with their common nucleus $F_0$ is a degree $\vert \C \vert +1$ maximal arc in $\PG(2,q)$.
\end{theorem}

Note that a maximal arc of degree $d$ of Mathon type contains Mathon sub-arcs of degree $d'$ for all $d'$ dividing $d$ (see \cite{MR1883870}).  As we mentioned above, Mathon's construction is a generalization of a previously known construction of Denniston.  This can be seen as follows.  Choose $\alpha \in \GF(q)$ such that $\Tr(\alpha)=1$.  Let $A$ be a subset of $\GF(q)^{\star}= \GF(q) \setminus \{0\}$ such that $A \cup \{0\}$ is closed under addition.  Then the point set of the conics 
\[
\K_A=\{F_{\alpha,1,\lambda}:\lambda \in A\}
\] 
together with the nucleus $F_0=(0,0,1)$ is the set of points of a degree $\vert A \vert +1$ maximal arc in $\PG(2,q)$.  This construction is exactly the definition of a maximal arc of Denniston type.  The conics in $\K_A$ are a subset of the \textit{standard pencil of conics} given by
\[
\{F_{\alpha,1,\lambda}:\lambda \in \GF(q)\}.
\]  
This pencil partitions the points of the plane, not on the line $z=0$ into $q-1$ disjoint  conics on the common nucleus $F_0=(0,0,1)$.  The line $z=0$ is often called the \textit{line at infinity} of the pencil and is denoted by $F_{\infty}$.  It has been proved by Mathon \cite{MR1883870} that a degree 4 Mathon arc is necessarily of Denniston type.  
However there are various families of Mathon maximal arcs known that are not of Denniston type.  Every Mathon arc that is not of Denniston type will be called a \emph{proper Mathon arc}. Actually, the most difficult part in checking that a given subset of conics of $\F$ is a maximal arc lies in checking whether the trace condition of Lemma \ref{discon} holds. In Section \ref{geom} we will present a more geometric approach to these arcs that allows us to overcome this problem. Furthermore, this geometric approach will be the key to our main result, which is the enumeration of the isomorphism classes of  Mathon $8$-arcs in $\PG(2,2^h)$, $h > 4$ and $h \neq 7$ prime.  The enumeration problem for Mathon arcs was first studied in \cite{MR2036477}, where bounds were derived for the number of isomorphism classes of Mathon arcs of ``big'' degree. The techniques of \cite{MR2036477} however failed for small degree arcs, and the enumeration of such arcs was left as an open problem.

It might be good to give at this stage an account of the known maximal arcs in Desarguesian projective planes of small order.

\subsubsection*{Maximal arcs in small Desarguesian planes}
\begin{enumerate}
\item The plane $\PG(2,8)$ has up to isomorphism only one maximal arc of degree $4$; it is  of Denniston type and is the dual of the regular hyperoval.
\item The plane $\PG(2,16)$ has up to isomorphism two maximal arcs of degree $8$: the dual of the regular hyperoval which is of Denniston type, and the dual of the Lunelli-Sce hyperoval which is of proper Mathon type. It has two isomorphism classes of maximal arcs of degree $4$, both of Denniston type and both self-dual.
\item The plane $\PG(2,32)$  has 6 isomorphism classes of hyperovals and hence the same number of  maximal arcs of degree $16$. As far as the maximal arcs of Denniston type are concerned, there is one of degree $4$, its dual of degree $8$, and the dual of the regular hyperoval which is a maximal arc of degree $16$.  Mathon gives in his original paper \cite{MR1883870} a construction of 3 maximal arcs of degree $8$ (and hence of 3 maximal arcs of degree $4$), which are not of Denniston type. In this paper we will prove  that there are no other maximal arcs of Mathon type  of degree $8$.
\end{enumerate}

\section{A geometric approach}\label{geom}

The following lemma was proved by Aguglia, Giuzzi and Korchmaros.  

\begin{lemma}[\cite{MR2443285}]  \label{aggiko}
Given any two disjoint conics $C_1$ and $C_2$ on a common nucleus.  Then there is a unique degree $4$ maximal arc of Denniston type containing $C_1 \cup C_2$.
\end{lemma}

We will generalize this to a synthetic version of Mathon's construction.

\begin{lemma}  \label{extline}
Given a degree $d < q/2$ maximal arc $M$ of Mathon type, consisting of $d-1$ conics on a common nucleus $n$, and a conic $C$ disjoint from $M$ with the same nucleus $n$,  there exists a line external to $M \cup C$.
\end{lemma}
\bew
 First we count the number of secants to $M$.  Since $(q+1)(q/d-1)+1$ is the number of external lines to $M$, the number of secants to $M$ is equal to
\begin{align*}
q^2+q+1-((q+1)\biggl(\frac{q}{d}-1\biggr)+1)=\biggl(\frac{d-1}{d}\biggr)q^2+\biggl(\frac{2d-1}{d}\biggr)q+1.
\end{align*}  
Next we count the number of lines that intersect both $M$ and $C$.  At first we will disregard the $q+1$ tangents to $C$, they will be added at the end.   Since the tangents to $C$ are disregarded, a secant line $l$ to both $C$ and $M$ must intersect $C$ in 2 points and $M$ in $d$ points.  This implies that the total number of secants to both $M$ and $C$ is equal to  
\[
\frac{1}{2}\biggl( \frac{(q+1)(d-1)+1}{d}-1\biggr)(q+1)+q+1=\biggl(\frac{d-1}{2d}\biggr)q^2+\biggl(\frac{3d-1}{2d}\biggr)q+1.
\]
We know that the number of lines intersecting $C$ is $(q+1)q/2+q+1$.  This means that the number of lines that intersect $C$ but do not intersect $M$ is 
\[
\frac{(q+1)q}{2}+q+1-\biggl(\biggl(\frac{d-1}{2d}\biggr)q^2+\biggl(\frac{3d-1}{2d}\biggr)q+1\biggr)=\frac{q^2}{2d}+\frac{q}{2d}.
\]
Finally we are able to count the number of secants to $M \cup C$.  We find
\begin{align*}
\biggl(\frac{d-1}{d}\biggr)q^2+\biggl(\frac{2d-1}{d}\biggr)q+1+\frac{q^2}{2d}+\frac{q}{2d}&=&\biggl(\frac{2d-1}{2d}\biggr)q^2+\biggl(\frac{4d-1}{2d}\biggr)q +1 \\
                                                                                                                                     &<&q^2+q+1.
\end{align*}
This proves that there exists an external line to $M \cup C$.
\evb

Using Lemma \ref{extline} we are able to prove the following result, which can be seen as a synthetic version of Mathon's construction.

\begin{theorem}[Synthetic version of Mathon's theorem]  \label{synthmathon}
Given a degree $d$ maximal arc $M$ of Mathon type, $d < q/2$, consisting of $d-1$ conics on a common nucleus $n$, and a conic $C_d$ disjoint from $M$ with the same nucleus $n$, then there is a unique degree $2d$ maximal arc of Mathon type containing $M \cup C_d$.
\end{theorem}
\bew
Denote the $d-1$ conics in the maximal arc $M$ by $C_1,C_2,C_3,...,C_{d-1}$.  Due to Lemma \ref{extline} we know there exists an external line $r$ to $M \cup C_d$.  We recoordinatize the plane $\PG(2,q)$ in such a way that the line $r$ now has equation $z=0$ and the common nucleus $n$ has coordinates $(0,0,1)$.  This provides us with the setting in which the conic $C_i$ has equation $\alpha_i x^2 +xy+\beta_i y^2 + \lambda_i z^2=0$.  \\ 
Next we define $\overline{C_i}:=\alpha_i \beta_i$.  It is clear that $\Tr(\overline{C_i})=1$, $\forall i=1,...,d$.  
We can now construct the degree $2d$ maximal arc containing $M \cup C_d$.  Let $C_i \oplus C_d:=C_{i+d}$ $\forall i=1,...,d-1$.  The construction used in the proof of Lemma \ref{aggiko}, which is based on Mathon, implies that $\Tr(\overline{C_{i+d}})=1$.  Due to Lemma \ref{discon} it follows that $C_i$, $C_d$ and $C_{i+d}$ are mutually disjoint.  \\
Next we need to check that the conics $C_i$ and $C_{j+d}$, $\forall i,j = 1,..., d-1$, are disjoint, i.e. $\Tr(\overline{C_i \oplus C_{j+d}})=1$.  Let $C_i \oplus C_j=C_k$, another conic which is defined in the closed set $M$, then
\begin{align*}
\Tr(\overline{C_i \oplus C_{j+d}})&= \Tr(\overline{C_i \oplus C_j \oplus C_d})  \\
						  &= \Tr(((\alpha_i \oplus \alpha_j) \oplus \alpha_d)((\beta_i \oplus \beta_j) \oplus \beta_d))  \\
						  &= \Tr((\alpha_k \oplus \alpha_d)(\beta_k \oplus \beta_d))  \\
						  &= \Tr(\overline{C_k \oplus C_d})  \\
						  &= \Tr(\overline{C_{k+d}})  \\
						  &= 1.	
\end{align*}       
Also the conics $C_{i+d}$, $\forall i=1,...,d-1$, have to be mutually disjoint.  This holds since
\begin{align*}
\Tr(\overline{C_{i+d} \oplus C_{j+d}})&= \Tr(\overline{C_i \oplus C_d \oplus C_j \oplus C_d})  \\
						  &= \Tr(\overline{C_i \oplus C_j})  \\
						  &= \Tr(\overline{C_{k}})  \\
						  &= 1,	
\end{align*}       
where again $C_k=C_i \oplus C_j$ is a conic in the original degree $d$ maximal arc $M$ of Mathon type.  It now follows that $\bigcup_{i=1}^{2d-1} C_{i}$ is a closed set of conics on a common nucleus $n$ which, due to Theorem \ref{mathon}, gives rise to a degree $2d$ maximal arc of Mathon type. It follows from Lemma  \ref{aggiko} and the above construction that this maximal arc is unique.
\evb

\section{Denniston $4$-arcs}

In \cite{MR1816247} Hamilton and Penttila determined the collineation stabiliser of a degree $d$ Denniston maximal arc.  

\begin{theorem}  \label{pentham}
In $\PG(2,2^e)$, $e > 2$, let $\D$ be a degree $d$ Denniston maximal arc, $q=2^e$, $2 < d < q/2$, with additive subgroup $A$.  Define the group $G$ acting on $\GF(2^e)$ by
\[
G=\{x \mapsto ax^{\sigma}: a \in \GF(2^e)^{\ast} , \sigma \in \Aut \> \GF(2^{2e})\}.
\]
Then the collineation stabiliser of $\D$ is isomorphic to $C_{2^e+1} \rtimes G_A$, the semidirect product of a cyclic group of order $(2^e+1)$ with the stabiliser of $A$ in $G$.
\end{theorem}

In the next lemma we will show that the order of $G_A$ is 2 in $\GF(2^{2h+1})$, $2h+1$ prime and $2h+1 \neq 2,3$, with $A$ the additive subgroup of a degree 4 Denniston maximal arc.

\begin{lemma}  \label{GA2}
In $\PG(2,2^{2h+1})$, $2h+1$ prime, and $2h+1 \neq 3$, let $\D$ be a degree $4$ Denniston maximal arc with additive subgroup $A$.  Define the group $G$ acting on $\GF(2^{2h+1})$ by
\[
G=\{x \mapsto ax^{\sigma}: a \in \GF(2^{2h+1})^{\ast} , \sigma \in \Aut \> \GF(2^{4h+2})\}.
\]
Then $\vert G_A \vert = 2$.
\end{lemma}
\bew
First we remark that the plane $\PG(2,2^{2h+1})$ can be coordinatized in such a way that the additive subgroup $A=\{0,1,w,w+1\}$, with $w \in \GF(2^{2h+1}) \setminus \{0,1\}$, is associated to the maximal arc $\D$ of Denniston type.  We will denote the multiplicative order of the element $w \in A$ in $\GF(2^{2h+1})$ by $o(w)$.  \\
Let $\varphi \in G_A$.  Since $\varphi (0)=0$ we can restrict the action of $\varphi$ on $A$ to its action on $\{1,w,w+1\}$.  The action of $\varphi$ on each element of $\{1,w,w+1\}$ has either order 1, 2 or 3.  \\
First we suppose $\sigma =1$. 

\begin{itemize}
\item If $a=1$ then $\varphi=\id$ in $G$.  
\item If $a \neq 1$ then the action of $\varphi$ on $1$ has either order 2 or 3. 
\begin{itemize}
\item If the order is 2 then 
\[
\varphi(\varphi(1))=a^2=1
\]  
which implies that $a=1$, clearly a contradiction.
\item If the order is 3 then
\[
\varphi(\varphi(\varphi(1)))=a^3=1
\]  
which implies that $3 \vert 2^{2h+1}-1$.  But since
\[
2^{2h+1}-1=2^{2h}+2^{2h-1}+ \cdots +1=2^{2h}+2^{2h-2}3+2^{2h-4}3+ \cdots +2^23+3,
\]  
we again  find a contradiction.
\end{itemize}
\end{itemize}
From now on suppose $\sigma \neq 1$.

\begin{enumerate}
\item[(1.)] Assume $\varphi$ acts trivially on $\{1,w,w+1\}$.  Then $\varphi(1)=1$ implies $a=1$.  Furthermore $\varphi(w)=aw^{\sigma}=w^{\sigma}$.  Since the action of $\varphi$ on each element of $\{1,w,w+1\}$ has order 1 there has to follow that $w^{\sigma}=w$, which implies $w^{\sigma-1}=1$.  This means $o(w) \vert \sigma -1$ but of course we know $o(w) \vert 2^{2h+1}-1$.  Now suppose $\sigma = 2^l$, $l \in \N^*$.  Note that $l < 4h+2$.  Then:
\[
o(w) \vert \ggd(2^l-1,2^{2h+1}-1),
\]     
which implies that
\[
o(w) \vert 2^{\ggd(l,2h+1)}-1.
\]
Now two possibilities can occur.
\begin{itemize}
\item If $l=2h+1$ then $\varphi: x \mapsto x^{2^{2h+1}}$, and so $\varphi$ indeed acts trivially on $A$.
\item If $l \neq 2h+1, 0$ then $\ggd(l,2h+1)=1$.  It follows that $o(w)=1$ and so $w=1$, which is clearly a contradiction.
\end{itemize} 

\item[(2.)] Assume the orbit on some element of $\{1,w,w+1\}$ has length 2 under the action of $\varphi$.  We consider two cases.
\begin{enumerate}
\item If $\varphi(1)=1$ then of course $a=1$ holds again.  This implies $\varphi(w)=w^{\sigma}$ and $\varphi(w^{\sigma})=w^{\sigma^2}$ but since the action of $\varphi$ has order 2 it follows that $w^{\sigma^2}=w$, implying $w^{\sigma^2-1}=1$.  We find that $o(w) \vert \sigma^2-1$ and also $o(w) \vert 2^{2h+1}-1$.  Using $\sigma = 2^l$ as we did above, we find, as $2\not\vert \>\> 2h+1$ and $2h+1$ is prime,
\[
o(w) \vert \ggd(2^{2l}-1,2^{2h+1}-1) \Rightarrow o(w) \vert 2^{\ggd(2l,2h+1)}-1 \Rightarrow o(w) \vert 2^{\ggd(l,2h+1)}-1.
\]   
Now the same two possibilities as in (1.) can occur, hence $\varphi$ acts trivially on $A$, clearly a contradiction.
\item Without loss of generality we can assume that $\varphi(1)=w$.  In this case we find that $a=w$.  Furthermore $\varphi(\varphi(1))=\varphi(w)=w^{\sigma +1}$ and so $w^{\sigma +1}=1$ since the action of $\varphi$ has order 2.  This implies $w^{\sigma^2-1}=1$ which gives us $o(w) \vert \sigma^2-1$ and again we know $o(w) \vert 2^{2h+1}-1$.  Using the same arguments as we did in (a), we see that 
\[
o(w)\vert  2^{\ggd(l,2h+1)}-1.
\]
Again the two possibilities we encountered in (1.) can occur.  
\begin{itemize}
\item If $l=2h+1$, then $\varphi: x \mapsto wx^{2^{2h+1}}$ and again
\[
\varphi(\varphi(1))=w^2=1,
\]  
a contradiction.
\item If $l \neq 2h+1, 0$ then $\ggd(l,2h+1)=1$.  It follows that $o(w)=1$ and so $w=1$, a contradiction.
\end{itemize}
\end{enumerate}

\item[(3.)] Now assume the orbit length is 3 under the action $\varphi$.  Without loss of generality we can assume that $\varphi(1)=w$, then $a=w$.  From this we find that $\varphi(\varphi(\varphi(1)))=w^{\sigma^2+\sigma+1}$, which of course has to be equal to 1.   We deduce that $w^{\sigma^3-1}=1$, implying that $o(w) \vert \sigma^3-1$ while $o(w) \vert 2^{2h+1}-1$ still holds.  If we again set $\sigma = 2^l$, $l \in \N^*$ and $l < 4h+2$, we find that $o(w) \vert 2^{\ggd(l,2h+1)}-1$, since $3 \!\!\not|\; 2h+1$.  Remark that in case $2h+1=3$ the degree 4 maximal arc would be a dual hyperoval of $\PG(2,8)$.  
The same two possibilities as in (1.) can occur.  \\
\begin{itemize}
\item If $l=2h+1$, then $\varphi: x \mapsto wx^{2^{2h+1}}$ and again
\[
\varphi(\varphi(\varphi(1)))=w^3=1,
\]  
a contradiction.
\item If $l \neq 2h+1, 0$ then $\ggd(l,2h+1)=1$.  It follows that $o(w)=1$ and so $w=1$, a contradiction. 
\end{itemize}
\end{enumerate}
We have proven that $\varphi$ either is $\id \in G$ or $\varphi: x \mapsto x^{2^{2h+1}}$, hence $\vert G_A \vert=2$.
\evb

\begin{remark}\label{aut4boog}
{\rm   We have just shown that if $q=2^p$, $p$ prime, $p\neq 2,3$, then the full automorphism group $G$ of a degree 4 Denniston arc has size $2(q+1)$ and is isomorphic to $C_{q+1}\rtimes C_2$. Let us have a closer look at the action of this group on the arc. It is well known that in $G$ there is a cyclic subgroup of order $q+1$ stabilizing all three conics of the arc and acting sharply transitively on the points of each of these conics. Furthermore this group stabilizes the line at infinity $L$ of the pencil determined by the arc and acts sharply transitively on the points of this line. The group $G$ also contains $q+1$ involutions. These involutions are exactly the $q+1$  elations with axis a line through the nucleus, and center the intersection of this line with the line at infinity $L$, stabilizing each of the three conics of the arc. There is exactly one such involution for each line through the nucleus. }
\end{remark}

In the following lemma we count the number of isomorphism classes of degree 4 maximal arcs of Denniston type.

\begin{lemma}
The number of isomorphism classes of degree $4$ maximal arcs of Denniston type  in $\PG(2,2^{2h+1})$, $2h+1$ prime, $2h+1\neq3$ is 
\[
N=\frac{2^{2h}-1}{3(2h+1)}.
\]
\end{lemma}
\bew
Since, by recoordinatizing the plane, we can always assume that a degree $4$ maximal arc of Denniston type is contained in the standard pencil, it suffices to calculate  the number of isomorphism classes of degree 4 maximal arcs in the standard pencil.

First of all we count the total number of degree 4 maximal arcs of Denniston type in the standard pencil.  We have $(2^{2h+1}-1)$ choices to pick a first conic and $(2^{2h+1}-2)$ choices to pick a second conic.  Since Lemma \ref{aggiko} states that there is a unique degree 4 maximal arc containing these 2 conics the total number of degree 4 maximal arcs in the standard pencil is
\[
\frac{(2^{2h+1}-1)(2^{2h+1}-2)}{6}. 
\]  
Let $\D$ be a degree 4 maximal arc of Denniston type.  Due to Theorem \ref{pentham} and Lemma \ref{GA2} we know that 
\[
\vert \Aut (\D) \vert = 2(2^{2h+1}+1).
\] 
Using this along with the fact that the order of the collineation stabiliser of the standard pencil is $2(2h+1)(2^{4h+2}-1)$ (see proof of Theorem \ref{pentham}), we can count the number of degree 4 maximal arcs of Denniston type that are isomorphic to $\D$.  We obtain
\[
\frac{2(2h+1)(2^{4h+2}-1)}{2(2^{2h+1}+1)}=(2h+1)(2^{2h+1}-1).
\]
Finally the number of isomorphism classes of degree 4 maximal arcs of Denniston type in the pencil is
\[
\frac{(2^{2h+1}-1)(2^{2^{2h+1}-2})}{6(2h+1)(2^{2h+1}-1)}=\frac{2^{2h}-1}{3(2h+1)}.
\]
\evb

\begin{lemma}  \label{isomden}
The number of degree {\rm4} maximal arcs of Denniston type in the standard pencil in $\PG(2,2^{2h+1})$, $2h+1$ prime, $2h+1\neq3$ which are isomorphic to a given one and contain a given conic $C$ equals $3(2h+1)$.
\end{lemma}
\bew
Let $\D$ be any degree $4$ maximal arc. The result follows immediately from the facts that the standard pencil contains $(2h+1)(2^{2h+1}-1)$ isomorphic copies of $\D$, the standard pencil contains $2^{2h+1}-1$ conics, and $\D$ contains $3$ conics, keeping in mind that $\mathrm{Aut}(\D)$ acts as described in Remark \ref{aut4boog}.
\evb

\section{Mathon $8$-arcs}

Let us first have a look at the geometric structure of a maximal $8$-arc of Mathon type; this is based on \cite{MR1940336}. Note that if $\K$ is a maximal arc constructed from a closed set of conics $\C$ on a common nucleus, then the point set of that arc contains no non-degenerate conics apart from those of $\C$ (see \cite{MR1997407}). From Lemma \ref{aggiko} it immediately follows that every Mathon $8$-arc contains exactly $7$ Denniston $4$-arcs, and each two of these seven $4$-arcs have exactly one conic in common.  One in fact easily sees that the structure with as point set the conics of $\K$, line set the degree $4$ subarcs of Denniston type, and the natural incidence is isomorphic to $\PG(2,2)$. In accordance with \cite{MR1940336} we define the lines at infinity of $\K$ to be the lines at infinity of each of the pencils determined by the degree $4$ subarcs. If $\K$ is of Denniston type there is a unique line at infinity, otherwise there are exactly $7$ distinct lines at infinity (see Theorem 2.2 of \cite{MR1940336} and the remark preceding it). Suppose namely that two subarcs $\K_1$ and $\K_2$ would have the same line at infinity. Let $C$ be the conic belonging to both $\K_1$ and $K_2$. Since a conic and a line uniquely determine a pencil, it follows that $\K_1$ and $\K_2$ belong to the same pencil, yielding that $\K$ is of Denniston type.  Note that it is essential here that any two of the degree $4$ arcs have a conic in common. 
In \cite{MR1940336} it is noticed that all known Mathon $8$-arcs seem to have an involution stabilizing $\K$ and all of its conics. Theorem 2.3 of \cite{MR1940336} gives a sufficient condition for such an involution to exist. In the next lemma we show that such an involution always exists.

\begin{lemma}  \label{iota}
Let $\K$ be a proper Mathon $8$-arc. Then the $7$ lines at infinity of $\K$ are concurrent and there exists a unique involution stabilizing $\K$ and all conics contained in $\K$. This involution is the elation with center the point of intersection of the lines at infinity and axis the line containing the nucleus of $\K$ and the center.
\end{lemma}
\bew
Denote the $7$ degree $4$ Denniston subarcs of $\K$ by $\D_i$, $i=1,\hdots,7$. Let $n$ be the nucleus of (the conics of) $\K$. Let $L_i$ be the line at infinity of $\D_i$. Let $c$ be the intersection of $L_1$ and $L_2$. Consider the unique involution $\iota$ with center $c$ and axis $nc$ that stabilizes the conic $C$ that is the intersection of $\D_1$ and $\D_2$. It is well known that $\iota$ will stabilize all conics in $\D_1$ and $\D_2$ (see e.g.\  the proof of Theorem \ref{pentham}). Now let $\D_3$ be the unique third $4$-arc that contains $C$. As $\K$ is uniquely determined by $\D_1$ and $\D_2$ (see Theorem \ref{geom}) it follows that $\iota$ must stabilize $\D_3$. Hence it must stabilize the line at infinity of $\D_3$, implying that $L_3$ contains $c$. It now also follows that $\iota$ stabilizes all conics of $\K$ and that all lines at infinity have to be stabilized;  we deduce that all lines at infinity are concurrent at $c$.
\evb 

\begin{corollary}  \label{C2}
Let $\K$ be a proper Mathon $8$-arc in $\PG(2,2^p)$, $p$ prime, $p\neq2,3,7$. Then $\mathrm{Aut}(\K)\cong C_2$.
\end{corollary}
\bew
Let $\phi$ be a non-trivial automorphism of $\K$. Clearly $\phi$ has to fix the intersection point $c$ of the lines at infinity of $\K$. 

First suppose that $\phi$ stabilizes one of the degree $4$ maximal subarcs of $\K$. From Remark \ref{aut4boog} and the fact that $c^\phi=c$ it follows that $\phi$ is the unique involution $\iota$ described in the previous lemma.

So, suppose that $\phi$ does not stabilize any of the Denniston subarcs. Hence no orbit of $\phi$ on the subarcs has length $1$. As there are $7$ subarcs, the set $O$ of orbit lengths has to be one of the following: $\{7\},\{5,2\},\{4,3\},\{3,2\}$. Suppose $O=\{3,2\}$. Then $\phi^2$ stabilizes some subarc and hence has to be the involution $\iota$. It follows that $\phi$ cannot have an orbit of length $3$, contradiction. The cases $O=\{5,2\}$ and $O=\{4,3\}$ are excluded in an analogous way. 

Hence $\phi$ cyclically permutes the $7$ subarcs. Suppose that $\phi$ would belong to $\mathrm{PGL}(3,2^p)$. As $\phi$ fixes the line $nc$ containing the nucleus and $c$, and $2^p$ is not divisible by $7$, we see that $\phi$ must fix a second line through $c$. If $\phi$ would fix a third line through $c$ it would fix all lines through $c$, a contradiction as $\phi$ cyclically permutes the lines at infinity of $\K$. Hence $7$ divides $2^p-1$, which implies that $3$ divides $p$, a contradiction. Hence $\phi\in\PGL(3,2^p)\setminus\mathrm{PGL}(3,2^p)$. As $7$ is prime it follows that $7$ divides the prime $p$, yielding that $p=7$, the final contradiction. 
\evb

In order to be able to count the number of isomorphism classes of degree 8 maximal arcs of Mathon type we need to know how many isomorphic images of a given degree 8 maximal Mathon arc there are. The following technical lemma will play a key role in our final calculations.

\begin{lemma}  \label{isom}
Let $\K$ be a proper Mathon $8$-arc  in $\PG(2,2^{2h+1})$, $2h+1$ prime, and $h\neq 1,3$. Then the number  of degree 8 maximal arcs  isomorphic to $\K$ that have one of their degree $4$ maximal subarcs in the standard pencil, contain a fixed given conic $C$ from the standard pencil and have the same intersection point for their lines at infinity is $21(2h+1)$.
\end{lemma}
\bew
Let $C$ be a conic in the standard pencil.  It is well known  that $G:=\Aut(C) \cong \PGL(2,2^{2h+1})$.  Hence $\left|G\right|=\vert \PGL(2,2^{2h+1}) \vert=(2h+1)(2^{2h+1}+1)(2^{4h+2}-2^{2h+1})$, which is the number of group elements that stabilize $C$ and its nucleus $n$.  The group $G$ acts transitively on the points not on $C$ and distinct from $n$. From this we can deduce that 
\[
\vert G_{C,n,(0,1,0)} \vert=\frac{(2h+1)(2^{2h+1}+1)(2^{4h+2}-2^{2h+1})}{(2^{2h+1}+1)(2^{2h+1}-1)}=(2h+1)2^{2h+1}.
\]    
The group $G_{C,n,(0,1,0)}$ acts transitively on the lines through $(0,1,0)$ that do not intersect $C$. Since $\frac{2^{2h+1}}{2}$ is the number of such lines, this implies that 
\[
\vert G_{C,[X=0],[Z=0]} \vert = \frac{\vert G_{C,n,(0,1,0)} \vert}{\frac{2^{2h+1}}{2}}=4h+2.
\]
Now suppose $\K$ is a proper Mathon arc of degree 8. Let $\D_i$, $i=1\hdots,7$ denote the seven $4$-arcs of Denniston type contained in $\K$, and let $C_1=C, \dots, C_7$ denote the seven conics of $\K$. Without loss of generality we may suppose that 
$\D_1$ belongs to the standard pencil and that $C$ is the conic belonging to both $\D_1,\D_2$ and $\D_3$. Furthermore we may assume that $(0,1,0)$ is the intersection point of the lines at infinity of $\K$. We want to count the number of isomorphic images of $\K$ that contain $C$, have a degree $4$ subarc in the standard pencil, and that have $(0,1,0)$ as intersection point of the lines at infinity.  Recall that $\left|\mathrm{Aut}(\K)\right|=2$. Let $\phi$ be an automorphism of the plane mapping $\K$ onto an isomorphic image of the desired type. First suppose $\phi$ stabilizes $C$ and the standard pencil. From the above we know that there are $4h+2$ choices for $\phi$.  Also, there are exactly $4h+2$ choices for $\phi$ that would map the pencil determined by $\D_i$, $i=2,3$, onto the standard pencil and stabilize $C$. We obtain $3(4h+2)$ choices for $\phi$ that stabilize $C$. 
Now let $C_i$, $i\neq1$ be any other conic of $\K$. Suppose that $C_i^\phi=C$. As one of the three pencils determined by $C_i$ and $\K$ has to be mapped onto the standard pencil, we see in an analogous way that there are $3(4h+2)$ choices for $\phi$ such that $C_i^\phi=C$. We obtain that in total there are $21(4h+2)$ choices for $\phi$. It follows that there are exactly $21(2h+1)$ isomorphic images of $\K$ of the desired type.    
\evb

Given a degree 4 maximal arc of Denniston type $\D_1$ in the standard pencil consisting of the conics $C_1, C_k, C_{k+1}$.  Due to Lemma \ref{aggiko} each conic $C$ disjoint from $\D_1$ together with $C_1$ will give rise to another degree 4 maximal arc of Denniston type which will be isomorphic to one of the degree 4 maximal arcs of Denniston type in the standard pencil.  In what follows we will establish the trace conditions that express the disjointness of the conic $C$ with respect to $\D_1$.  \\
Let $\D_1$ and $\D_2$ be 2 non-isomorphic degree 4 maximal arcs of Denniston type.  Without loss of generality we can assume that both arcs are contained in the standard pencil and that both contain a common conic $C_1$.    
Let the additive subgroups $\{0,1,k,k+1\}$ and $\{0,1,l,l+1\}$, with $k \neq l,l+1$ and $k,l \in \GF(2^{2h+1}) \setminus \{0,1\}$, be the ones associated to the maximal arcs $\D_1$ and $\D_2$ respectively.  In other words we assume $\D_1$ consists of the conics $C_i, i=1,k,k+1$ given by the equation
\[
C_i :x^2+xy+y^2+i z^2=0
\]
and $\D_2$ consists of the conics $C_j, j=1,l,l+1$ given by
\[
C_j :x^2+xy+y^2+j z^2=0.
\]
Consider the automorphisms  $\theta$ of $\PG(2,2^{2h+1})$ determined by the matrix

\begin{equation}  \label{autom}
	\begin{pmatrix}
		\sqrt{\lambda}^{-\sigma}&0&0\\t&\sqrt{\lambda}^{-\sigma}&0\\ \sqrt{\sqrt{\lambda}^{-\sigma}t+t^2}&0&1
	\end{pmatrix},
\end{equation}

\noindent and the field automorphism $\sigma$, with $\lambda=1,l,l+1$ and $t \in \GF(2^{2h+1})$.  These automorphisms will map $C_\lambda$ onto $C_1$ while $(0,0,1)^\theta=(0,0,1)$ and $(0,1,0)^\theta =(0,1,0)$.  In fact all automorphisms of $\PG(2,2^{2h+1})$ which fix $(0,0,1)$ and $(0,1,0)$ and map $C_{\lambda}$ onto $C_1$ are of the form $\theta$.  There are three possibilities for $\theta$ that we have to take into account: $C_1^{\theta}=C_1$, $C_l^{\theta}=C_1$ and $C_{l+1}^{\theta}=C_1$. We will look at the case where $C_l$ is mapped onto $C_1$ and examine what values for $t$ satisfy the conditions
\[
C_1^{\theta} \cap C_k =\emptyset
\] 
and
\[
C_1^{\theta} \cap C_{k+1} =\emptyset.
\]
Analogous results can be found in the cases $C_1^{\theta}=C_1$ and $C_{l+1}^{\theta}=C_1$.  First we construct the image of $C_1$ under $\theta$.  It is clear that the point $(0,1,1)$, which is the intersection of $C_1$ and the $x$-axis, is mapped onto the point $(0,\sqrt{l}^{-\sigma},1)$.  It is clear that $(0,\sqrt{l}^{-\sigma},1) \neq (0,\sqrt{k},1)$, since $l^{-\sigma} = k$  would immediately imply that $C_1^{\theta} \cap C_k \neq \emptyset$, a contradiction.  Analogously $(0,\sqrt{l}^{-\sigma},1) \neq (0,\sqrt{k+1},1)$, i.e. $l^{-\sigma} \neq k+1$, since in this case the contradiction $C_1^{\theta} \cap C_{k+1} \neq \emptyset$ would hold.  \\
 
Furthermore we look at the image of a general point $(1,y,z)$ of $C_1$, $y,z \in \GF(2^{2h+1})$, where of course $1+y+y^2+z^2=0$ holds.  We find 
\[
	\begin{pmatrix} 
		\sqrt{l}^{-\sigma}&0&0\\t&\sqrt{l}^{-\sigma}&0\\ \sqrt{\sqrt{l}^{-\sigma}t+t^2}&0&1
	\end{pmatrix}
	\begin{pmatrix}
		1\\y\\z
	\end{pmatrix}^{\sigma}=
	\begin{pmatrix}
		\sqrt{l}^{-\sigma}\\t+\sqrt{l}^{-\sigma}y^{\sigma}\\ \sqrt{\sqrt{l}^{-\sigma}t+t^2}+z^{\sigma}
	\end{pmatrix},
\]
with $\sigma \in \Aut(\GF(2^{2h+1}))$.  The condition $C_1^{\theta} \cap C_k =\emptyset$ is satisfied if and only if the equation
\[
l^{-\sigma}+\sqrt{l}^{-\sigma}t+l^{-\sigma}y^{\sigma}+t^2+l^{-\sigma}y^{2\sigma}+k\sqrt{l}^{-\sigma}t+kt^2+kz^{2\sigma}=0
\] 
has no solutions in $\GF(2^{2h+1})$.  Equivalently, since $1+y^{\sigma}+y^{2\sigma}=z^{2\sigma}$, we find
\begin{align*}
(l^{-\sigma}+k)z^{2\sigma}+\sqrt{l}^{-\sigma}t+t^2+k\sqrt{l}^{-\sigma}t+kt^2&=&0  \\
\Leftrightarrow z^{2\sigma}&=&\frac{(1+k)t(\sqrt{l}^{-\sigma}+t)}{(l^{-\sigma}+k)}.
\end{align*}
Hence the conics $C_1^{\theta}$ and $C_k$ will be disjoint if and only if the equation
\[
1+y^{\sigma}+(y^{\sigma})^2+\frac{(1+k)t(\sqrt{l}^{-\sigma}+t)}{(l^{-\sigma}+k)}=0.
\]
has no solutions in $y^{\sigma}$, or equivalently if and only if 
\[
\Tr\Big[1+\frac{(1+k)t(\sqrt{l}^{-\sigma}+t)}{(l^{-\sigma}+k)}\Big]=1.
\]
Since $\Tr(1)=1$ in $\GF(2^{2h+1})$ we find the condition

\begin{align}  \label{trace1}
\Tr\Big[\frac{(1+k)t(\sqrt{l}^{-\sigma}+t)}{(l^{-\sigma}+k)}\Big]=0.
\end{align}
Analogously, the trace condition 

\begin{align}  \label{trace2}
\Tr\Big[\frac{kt(\sqrt{l}^{-\sigma}+t)}{(l^{-\sigma}+k+1)}\Big]=0
\end{align}

is necessary and sufficient for $C_1^{\theta} \cap C_{k+1}=\emptyset$.  \\
It is clear that also the conic $C_{l+1}^{\theta}$ has to be disjoint from both $C_k$ and $C_{k+1}$.  However, due to Lemma \ref{aggiko}, we know that the two conics $C_1$ and $C_1^{\theta}$ give rise to a unique degree 4 maximal arc of Denniston type.  The third conic contained in this 4-arc has to be $C_{l+1}^{\theta}$, since we are actually looking at the image of $\D_2$ under $\theta$.  Using Theorem \ref{synthmathon} we know that both the degree 4 maximal arcs $\D_1$ and the conic $C_1^{\theta}$ induce a unique degree 8 maximal arc in which of course all conics are mutually disjoint.  Since $\D_2^{\theta}$ is contained in this 8-arc we can conclude that $C_{l+1}^{\theta}$ will be disjoint from all other conics in the 8-arc.  This implies that the two trace conditions originating from the disjointness of $C_{l+1}^{\theta}$ will lead to the same values for $t$.  

Next, consider a degree 4 maximal arc $\D$ in the degree 8 maximal arc.  If $\theta_{t'}=\iota \theta_t$, where $\iota$ is the unique involution described in Lemma \ref{iota}, fixing all conics in the 8-arc, then we know $\D^{\theta_t}=\D^{\theta_{t'}}$.  Since $\theta_{t'} \neq \theta_t$, the values $t$ and $t'$ will of course be distinct.  However, these $t$-values have to give rise to the same degree 4 arc $\D^{\theta_t}$.  In other words, these $t$-values come in pairs, which means that two $t$-values induce one and the same line at infinity or equivalently, one and the same degree 4 maximal arc of Denniston type.  

Suppose there would be a third value $t''$ inducing the same degree 4 arc of Denniston type.  This means $\D^{\theta_t}=\D^{\theta_{t''}}$ or $\D^{\theta_t \theta_{t''}^{-1}}=\D$.  Since $t$ and $t''$ are presumed to be distinct, it follows that $\theta_t \theta_{t''}^{-1}= \iota$ which means that $\theta_t=\iota \theta_{t''}$ or equivalently $\iota \theta_t=\theta_{t''}$.  We conclude that $\theta_{t''}=\theta_{t'}$ or $t''=t'$.

\begin{remark}  \label{identsigma}
{\rm
There are no restrictions on $\sigma$ since $\D_1$ and $\D_2$ are non-isomorphic.  On the other hand, consider $\D_1$ consisting of the conics $C_1,C_k,C_{k+1}$ and the automorphism fixing the conic $C_1$.  If in that case $\sigma$ is the identity then the conics $C_k$ and $C_k^{\theta}$ will intersect in the point $(0,\sqrt{k},1)$ on the $x$-axis.  Analogously the conics $C_{k+1}$ and $C_{k+1}^{\theta}$ intersect in $(0,\sqrt{k+1},1)$.  Of course this does not occur in disjoint conics. 
}
\end{remark}
	
\begin{theorem}  \label{8count}
The number of isomorphism classes of degree 8 proper Mathon arcs in $\PG(2,2^{2h+1})$, $2h+1 \neq 7$ and prime, is exactly
\[
\frac{N}{14}(2^{2h-2}-1)((6h+3)N-1),
\]
where $N=\displaystyle{\frac{(2^{2h}-1)}{3(2h+1)}}$.
\end{theorem}
\bew
Let $\D^i$, $i=1,...,N$, be chosen fixed and representative of each isomorphism class of degree 4 maximal arcs of Denniston type in the standard pencil.  Assume $\D^i$ consists of the conics $C_1,C_2^i$ and $C_3^i$, $i=1,...,N$.  First of all we want to calculate how many degree 8 maximal arcs of Mathon type contain one of the $N$ degree 4 maximal arcs $\D^i$, say $\D^1$, have the $x$-axis as elation axis and the intersection point of the lines at infinity as elation centre.  \\
Assume $i \neq 1$.  \\
Let $\theta$ be an automorphism  of $\PG(2,2^{2h+1})$ as given by the matrix  in (\ref{autom}).  We need to count in how many ways we can map $C_2^i$ onto $C_1$ such that both conditions
\begin{align*}
\left\{
	\begin{array}{l} 
		C_1^{\theta} \cap C_2^1=\emptyset\\C_1^{\theta} \cap C_3^1=\emptyset
	\end{array}
\right.
\end{align*}
are satisfied.  As seen above these conditions of disjointness are equivalent to the two trace conditions
\begin{align*}
\left\{
	\begin{array}{l} 
		\Tr[A_1(\sigma)t+B_1(\sigma)t^2]=0\\ \Tr[A_2(\sigma)t+B_2(\sigma)t^2]=0,
	\end{array}
\right.
\end{align*}
where $A_1,A_2,B_1$ and $B_2$ are functions of $\sigma$.  This can also be written as 
\begin{align*}
\left\{
	\begin{array}{l} 
		\Tr[(A_1(\sigma)+\sqrt{B_1(\sigma)})t]=0\\ \Tr[(A_2(\sigma)+\sqrt{B_2(\sigma)})t]=0,
		\end{array}
\right.
\end{align*}
which are two linear equations that correspond to two hyperplanes in the vector space $\V(2h+1,2)$.  Since $A_1(\sigma)+\sqrt{B_1(\sigma)} \neq A_2(\sigma)+\sqrt{B_2(\sigma)}$, which is easily checked by adding (\ref{trace1}) and (\ref{trace2}), the corresponding hyperplanes intersect in a $(2h-1)$-dimensional subspace.  We conclude that there are $2^{2h-1}=\frac{2^{2h+1}}{4}$ solutions to the system of trace conditions above.  This means that for every $\sigma$ there are $\frac{2^{2h+1}}{4}$ solutions for $t$.  However, since these $t$-values come in pairs we find, for every field automorphism $\sigma$, that there are $\frac{2^{2h+1}}{8}$ degree 4 maximal arcs.  One of them will give rise to a degree 8 maximal arc of Denniston type and so there are 
\[
(2h+1)\Big(\frac{2^{2h+1}}{8}-1\Big)
\]     
automorphisms $\theta$ that satisfy the needed conditions and induce a degree 8 maximal arc of Mathon type.  One such automorphism leads to two conics disjoint from $C_2^1$ and $C_3^1$ and so we get 
\[
(2h+1)\Big(\frac{2^{2h+1}}{4}-2\Big)
\]
conics disjoint from $C_2^1$ and $C_3^1$.  \\
In exactly the same way we can map $C_3^i$ onto $C_1$ and also $C_1$ onto $C_1$.  This gives us
\[
3(2h+1)\Big(\frac{2^{2h+1}}{4}-2\Big)
\]
conics that expand $\D^1$ to a degree 8 maximal arc of Mathon type.  \\
Now assume $i=1$.  \\
In the cases where $C_2^1$ is mapped onto $C_1$ and $C_3^1$ is mapped onto $C_1$ we  find again
\[
(2h+1)\Big(\frac{2^{2h+1}}{4}-2\Big)
\]
conics to expand $\D^1$.  If we consider the case where $C_1$ is fixed however, we have to make sure that $\sigma$ is not the identity as seen in the remark above.  And so in the case $i=1$ we get
\[
2(2h+1)\Big(\frac{2^{2h+1}}{4}-2\Big)+2h\Big(\frac{2^{2h+1}}{4}-2\Big)
\]
conics to expand $\D^1$.  As there are $N-1$ choices for $\D^i$, $i \neq 1$ there are a total of 
\[
(N-1)(6h+3)\Big(\frac{2^{2h+1}}{4}-2\Big)+(6h+2)\Big(\frac{2^{2h+1}}{4}-2\Big)
\]
such conics.  Suppose we counted one of these conics, say $C$, twice.  Since, due to Lemma \ref{aggiko}, this conic $C$ induces a unique degree 4 maximal arc together with $C_1$ it would imply that $C$ is the image of two conics contained in one of the $N$ 4-arcs $\D^i$.  However, this would give rise to an automorphism of the 4-arc that does not fix the conics, clearly a contradiction.  

In other words, we can use each one of these conics to expand $\D^1$ to a degree 8 maximal arc of Mathon type.  
Moreover, since the four conics disjoint from $\D^1$ in a degree 8 maximal arc of Mathon type all give rise to this same degree 8 arc, we find
\[
\frac{1}{4}\Big[(N-1)(6h+3)\Big(\frac{2^{2h+1}}{4}-2\Big)+(6h+2)\Big(\frac{2^{2h+1}}{4}-2\Big)\Big]
\] 
degree 8 maximal arcs of Mathon type that contain $\D^1$.  Of course there were $N$ choices for $\D^1$ and so there are
\[
\frac{N}{4}\Big[(N-1)(6h+3)\Big(\frac{2^{2h+1}}{4}-2\Big)+(6h+2)\Big(\frac{2^{2h+1}}{4}-2\Big)\Big]
\]
degree 8 maximal arcs of Mathon type that contain the degree 4 maximal arc $\D^i$.  As a result of Lemma \ref{isom} we now find 
\[
\frac{N}{28}\Big[(N-1)(6h+3)\Big(\frac{2^{2h+1}}{4}-2\Big)+(6h+2)\Big(\frac{2^{2h+1}}{4}-2\Big)\Big]
\]
isomorphism classes of degree 8 maximal arcs of Mathon type in $\PG(2,2^{2h+1})$, $2h+1 \neq 7$.  Remark that we divided by 7 as Lemma \ref{isomden} and Lemma \ref{isom} state.  This is due to the fact that we now fix an entire degree 4 maximal arc in the pencil, not only the conic $C_1$.   
\evb

\begin{remark}
{\rm
If $2h+1=7$ the situation changes.  Let $\phi$ be a non-trivial automorphism of $\K$.  

If $\phi$ stabilizes one of the degree 4 maximal subarcs of $\K$ we have seen in the proof of Corollary \ref{C2} that $\phi$ must be the unique involution $\iota$ described in Lemma \ref{iota}.

Now suppose that $\phi$ does not stabilize any of the Denniston subarcs.  Since 7 is the only possible orbit length of $\phi$ on these subarcs it follows that the order of $\langle \phi \rangle$ has to be a multiple of 7.  Let the order of $\langle \phi \rangle$ be $k7$, with $k \in \N^{\star}$.  In that case $\vert \langle \phi \rangle_{\D} \vert=k$.  Furthermore, since $\vert \Aut(\D) \vert =2$ we find that $k=2$.  This means that $\vert \Aut(\K) \vert=14$ and we can no longer benefit from the fact that $\Aut(\K) \cong C_2$, which implies that the previous counting arguments no longer hold.       
}
\end{remark}

\section{Maximal arcs in $\PG(2,32)$}

In this section we will consider the case $\PG(2,32)$.  Due to a randomized computer search Mathon found three isomorphism classes of  degree 8 maximal arcs in $\GF(32)$.  It now follows from Theorem \ref{8count} that there are exactly 3 such arcs.  In this section we will compose these arcs and conclude with the actual equations of their conics as they were written down by Mathon in \cite{MR1883870}.  In \cite{MR1816247} Hamilton and Penttila showed that there is a unique degree 4 maximal arc of Denniston type in $\PG(2,32)$, up to isomorphism.  Let $w$ be a primitive element in $\GF(32)$ satisfying $w^{18}+w=1$.  The three conics $C_1, C_w$ and $C_{w+1}$, given by
\[
\{x^2+xy+y^2+\lambda z^2 \vert \lambda \in \langle 1,w \rangle \setminus \{0\}\},
\]         
determine a degree 4 maximal arc of Denniston type $\D_1$ on the nucleus $(0,0,1)$.  \\
Due to the above the number of isomorphism classes of degree 8 maximal arcs of Mathon type in $\PG(2,32)$ is equal to the number of isomorphism classes of degree 8 maximal arcs of Mathon type that contain $\D_1$ while the intersection point $(0,1,0)$ of the lines at infinity is fixed.  This means we need to count the number of conics with nucleus $(0,0,1)$ that are disjoint from $\D_1$ while fixing the point $(0,1,0)$.  It is clear (Lemma \ref{aggiko}) that every such conic, together with the conic $C_1$, determines a degree 4 maximal arc of Denniston type $\D_2$, which of course is isomorphic to $\D_1$.  We now consider automorphisms $\theta$ of $\PG(2,32)$ such that $(\D_1)^{\theta}$ contains $C_1$, $(0,0,1)^{\theta}=(0,0,1)$ and $(0,1,0)^{\theta}=(0,1,0)$. We need to take in account three possibilities for $\theta$, more precisely: $C_1^{\theta}=C_1$, $C_w^{\theta}=C_1$ and $C_{w+1}^{\theta}=C_1$.  First let us consider the automorphism $\theta$ given by 
\[
\theta: x \rightarrow Mx^{\sigma},
\]
with 
\begin{equation*}
M:=\begin{pmatrix} 
		w^{-9\sigma}&0&0\\t&w^{-9\sigma}&0\\ \sqrt{w^{-9\sigma}t+t^2}&0&1
	\end{pmatrix},
\end{equation*} 
where $\sigma \in \Aut(\GF(32))$ and $t \in \GF(32)$.  This automorphism will indeed fix the points $(0,0,1)$ and $(0,1,0)$ while $C_{w+1}^{\theta}=C_1$.  The trace conditions that satisfy the conditions of disjointness $C_1^{\theta} \cap C_w= \emptyset$ and $C_1^{\theta} \cap C_{w+1}= \emptyset$ are

\[
\begin{cases} \displaystyle \Tr\Big[\frac{w^{9\sigma}t(1+w)(1+w^{9\sigma}t)}{(1+w^{1+18\sigma})}v\Big]=0 \\[2mm]  \displaystyle \Tr\Big[\frac{w^{9\sigma}t(1+w^{18})(1+w^{9\sigma}t)}{(1+w^{18+18\sigma})}\Big]=0. 
\end{cases}
\]

For all $\sigma \in \Aut(\GF(32))$ we find 8 elements $t \in \GF(32)$ satisfying these conditions.  More precisely, for every $\sigma$, we find the following $t$-values.  \\
\[
\begin{tabular}{lll} 
$\sigma=1:$ & $t=0,w^8,w^{22},w^{21},w^{11},w^{30},w^6,w^{15}$  \\ 
$\sigma=2:$ & $t=0,w^{13},w^{6},w^{28},w^{29},w^{22},w^{18},w^{15}$ \\
$\sigma=4:$ & $t=0,w^2,w,w^{19},w^{10},w^{22},w^{17},w^{26}$  \\
$\sigma=8:$ & $t=0,w^{21},w^{2},w^{13},w^{18},w^{16},w^{11},w^{15}$  \\
$\sigma=16:$ & $t=0,w^7,w^{9},w^{12},w^{29},w^{14},w^{17},w^{11}$
\end{tabular} 
\]
\smallskip

These 8 elements $t$ are partitioned into pairs.  For example if $\sigma =1$ we find the pairs
\begin{align}  \label{pairs}
(0,w^{22}),(w^8,w^{21}),(w^{11},w^{30}),(w^6,w^{15}).
\end{align}

The case $C_w^{\theta}=C_1$ can be handled in an analogous way.  The trace conditions now are

\[
\begin{cases} \displaystyle \Tr\Big[\frac{w^{-15\sigma}t(1+w)(1+w^{-15\sigma}t)}{(1+w^{\sigma+1})}\Big]=0 \\[2mm]  \displaystyle \Tr\Big[\frac{w^{-15\sigma}t(1+w^{18})(1+w^{-15\sigma}t)}{(1+w^{\sigma+18})}\Big]=0. 
\end{cases}
\]

The $t$-values for every $\sigma$, which are again partitioned in pairs, are listed below.
\[
\begin{tabular}{lll} 
$\sigma=1:$ & $t=0,w^{21},w^{19},w^{24},1,w^{25},w^{11},w^{15}$  \\ 
$\sigma=2:$ & $t=0,w^{20},w^{30},w^{24},w^{10},w^{14},w^{18},w^{23}$ \\
$\sigma=4:$ & $t=0,w^3,w^2,w^{20},1,w^{29},w^{5},w^{8}$  \\
$\sigma=8:$ & $t=0,w^{4},w^{12},w^{24},w^{5},w^{22},w^{27},w^{16}$  \\
$\sigma=16:$ & $t=0,w^4,w^{6},w^{9},w^{5},w^{22},w^{23},w^{15}$
\end{tabular} 
\]
\smallskip

Finally we have a look at the case where $C_1$ is fixed.  In accordance to Remark \ref{identsigma} we must demand that $\sigma \neq 1$ otherwise the $x$-axis is fixed pointwise and it would be impossible for the conics $C_w$ and $C_{w+1}$ to obtain disjoint images.  In the same way as seen above the conditions of disjointness result in the following system of trace conditions:

\[
\begin{cases} \displaystyle \Tr\Big[\frac{t(1+w)(1+t)}{(w+w^{\sigma})}\Big]=0 \\[2mm]  \displaystyle \Tr\Big[\frac{t(1+w^{18})(1+t)}{(w^{18}+w^{\sigma})}\Big]=0. 
\end{cases}
\]

The $t$-values for every $\sigma$ are:
\[
\begin{tabular}{lll} 
$\sigma=2:$ & $t=0,w^{7},w^{6},w^{24},1,w^{22},w^{27},w^{15}$ \\
$\sigma=4:$ & $t=0,w^4,w^{12},w^{24},1,w^{10},w^{23},w^{15}$  \\
$\sigma=8:$ & $t=0,w^{2},w,w^{19},1,w^{5},w^{18},w^{11}$  \\
$\sigma=16:$ & $t=0,w^4,w^{12},w^{24},1,w^{10},w^{23},w^{15}$
\end{tabular} 
\]
\smallskip

Each one of these pairs (e.g. \ref{pairs}) give rise to a unique degree 4 maximal arc of Denniston type.  This means that, for each one of them, we get two conics disjoint from $\D_1$.  One of these degree 4 maximal arcs is contained in the pencil of $\D_1$ and so it leads to a degree 8 maximal arc of Denniston type.  The other three induce proper Mathon arcs of degree 8.  \\
Now we are able to count the conics that give rise to a maximal arc of Denniston type (``$D$-conics") as well as the conics that give rise to a maximal arc of Mathon type (``$M$-conics").  Remark that only 4 values for $\sigma$ can be included in the case where $C_1$ is fixed since the identity leads to a contradiction.\\
\[
\begin{tabular}{lcc} 
& \multicolumn{1}{c}{``D-conics"} &\multicolumn{1}{c}{``M-conics"}\\ 
$C_{w+1}^{\theta}= C_1$ & $5 \times 1\times 2$ & $5 \times 3 \times 2$ \\ 
$C_{w}^{\theta}= C_1$ & $5 \times 1 \times 2$ & $5 \times 3 \times 2$ \\
$C_1^{\theta}= C_1$ & $4 \times 1\times 2$ & $4 \times 3 \times 2$ \\
& 28 & 84
\end{tabular} 
\]   
It follows that we find 84 ``M-conics".  This means there are 21 degree 8 maximal arcs of Mathon type.  As each of these proper Mathon arcs of degree 8 have an automorphism group of size 2 and contain exactly 7 degree 4 maximal arcs of Denniston type, which are isomorphic to $\D_1$, we obtain three isomorphism classes of  proper Mathon arcs of degree 8.  \par

On a more technical note we can calculate the equation of the conic $C_1^{\theta}$ using the matrix $M$ and the matrix
\begin{equation*}
A:=\begin{pmatrix} 
		1&1&0\\0&1&0\\0&0&1
	\end{pmatrix},
\end{equation*} 
associated to the equation $x^2+xy+y^2+z^2=0$ of $C_1$.  Analogous results hold for the conics $C_w^{\theta}$ and $C_{w+1}^{\theta}$.  This will enable us to construct the entire degree 8 maximal arc using Theorem \ref{synthmathon}.  We need to calculate the form $M^{T^{-1}}A^{\sigma}M^{-1}$.  Since $A=A^{\sigma}$ and 
\begin{equation*}
M^{-1}:=\begin{pmatrix} 
		w^{9\sigma}&0&0\\tw^{18\sigma}&w^{9\sigma}&0\\ \sqrt{w^{-9\sigma}t+t^2}w^{9\sigma}&0&1
	\end{pmatrix}
\end{equation*} 
we find that $(M^{-1})^TA^{\sigma}M^{-1}$ is equal to the matrix
\begin{equation*}
	\begin{pmatrix} 
		w^{18\sigma}+tw^{18\sigma}(w^{9\sigma}+tw^{18\sigma})+(w^{-9\sigma}t+t^2)w^{18\sigma}&w^{18\sigma}+tw^{27\sigma}&\sqrt{w^{-9\sigma}t+t^2}w^{9\sigma}\\tw^{27\sigma}&w^{18\sigma}&0\\\sqrt{w^{-9\sigma}t+t^2}w^{9\sigma}&0&1
	\end{pmatrix}.
\end{equation*} 
This means that the equation of the conic $C_1^{\theta}$ is given by
\[
(w^{18\sigma}+tw^{27\sigma}+t^2w^{36\sigma}+(w^{-9\sigma}t+t^2)w^{18\sigma})x^2+w^{18\sigma}xy+w^{18\sigma}y^2+z^2=0,
\]
with $t \in \GF(32)$ and $\sigma \in \Aut(\GF(32))$, which is equivalent to the equation
\[
(1+(1+w^{18\sigma})w^{-9\sigma}t+(1+w^{18\sigma})t^2)x^2+xy+y^2+w^{13\sigma}z^2=0.
\]
Let us now consider the case $\sigma=4$ and $t=w^{2}$.  We obtain 
\[
w^{12}x^2+xy+y^2+w^{21}z^2=0
\]
as the equation of $C_1^{\theta}$.  If we multiply this equation by $w^{19}$, set $y=w^{12}y'$ and $z=w^8z'$, we find 
\[
x^2+xy'+w^{12}y'^2+w^{25}z'^2=0,
\]
which is equivalent to 
\[
x^2+xy+w^{12}y^2+w^{25}z^2=0.
\]
Using Theorem \ref{synthmathon} and Mathon's composition we can easily compose the remaining three conics of the degree 8 arc.  Their equations are 
\begin{align*}
C_1 \oplus C_1^{\theta}&:& x^2+xy+w^{6}y^2+w^{21}z^2=0  \\
C_w \oplus C_1^{\theta}&:& x^2+xy+w^{18}y^2+w^{16}z^2=0  \\
C_{w+1} \oplus C_1^{\theta}&:& x^2+xy+w^{20}y^2+w^{9}z^2=0.
\end{align*}
This way we managed to construct the degree 8 maximal arc consisting of the conics $\{C_1,C_w,C_{w+1},C_1^{\theta},C_1 \oplus C_1^{\theta},C_w \oplus C_1^{\theta},C_{w+1} \oplus C_1^{\theta}\}$.  In \cite{MR1883870} Mathon found the three degree 8 maximal arcs (not of Denniston type) in $\PG(2,32)$ formed by 
\[
\{x^2+xy+(w^k+w^l\lambda+w^m\lambda^3)y^2+\lambda z^2 \vert \lambda \in \langle1,w,w^9 \rangle \setminus \{0\}\},
\]
with exponents $(k,l,m)=(12,15,4),(5,25,14)$, and $(6,19,8)$, respectively.  The 8-arc constructed above is exactly the one of Mathon corresponding to the exponents $(k,l,m)=(6,19,8)$.  The other two proper Mathon 8-arcs in $\GF(32)$ are found in an analogous way.  

\section*{Acknowledgement} 
Part of this paper has been written while the authors were visiting the Department of Mathematics of the University of California, San Diego.  We would also like to mention that this research was partly financed by the Research Project  ``Incidence Geometry''  (BOF/GOA/010) at Ghent University. 

\bibliographystyle{plain}

\end{document}